\documentclass[a4paper, 12pt]{amsproc}
\usepackage{amssymb}
\title[Hecke Algebra of Type $A_n$]
{A New Proof for Classification of Irreducible Modules
 \\  of a Hecke Algebra of Type $ A_n$}

\author[N. Xi]{Nanhua XI$^{*}$}
\address{$^{*}$
Institute of Mathematics\\
Chinese Academy of Sciences\\
Beijing, 100080\\
China } \email{nanhua@math.ac.cn}
\thanks{N. Xi was partially supported by Natural Sciences Foundation
of China (No. 10671193).} \dedicatory{Dedicated to Professor Gus
Lehrer on his sixtieth birthday}

\begin{document}
\baselineskip=18pt
\begin{abstract}
In this paper we  give a new proof for the classification of
irreducible modules of an affine Hecke algebra of type $A_n$, which
was obtained by G. E. Murphy in 1995.
\end{abstract}

\maketitle

\def\Cal{\mathcal}
\def\bold{\mathbf}
\def\ca{\mathcal A}
\def\cdz{\mathcal D_0}
\def\cd{\mathcal D}
\def\cdo{\mathcal D_1}
\def\bold{\mathbf}
\def\l{\lambda}
\def\le{\leq}
\def\o{\omega}

Let $H$ be the Hecke algebra of the symmetric group $S_n$ over a
commutative ring $K$ with an invertible parameter $q\in K$. In [M]
Murphy worked out a classification of irreducible modules of $H$
when $K$ is a field, which is similar to the classification of
irreducible modules of a symmetric group over a field [J]. In this
paper we shall give a new proof for Murphy's classification.
Essentially the idea is due to Murphy, but we use Kazhdan-Lusztig
theory and affine Hecke algebra of type $\tilde A_{n-1}$ to prove
his result by a direct calculation.

As usual, the simple reflections of $S_n$ consisting of the
transposes $s_i=(i,i+1)$ for $i=1,2,...,n-1$. As a free $K$-module,
the Hecke algebra $H$ has a basis $T_w,\ w\in S_n$, and the
multiplication is  defined by the relations $(T_s-q)(T_s+1)=0$ if
$s$ is a simple reflection, $T_wT_u=T_{wu}$ if $l(wu)=l(w)+l(u)$,
here $l:\ S_n\to \bold N$ is the length function.

For each partition $\lambda=(\l_1,...,\l_k)$ of $n$, set
$I_j=\{\l_1+\cdots+\l_{j-1}+1,\l_1+\cdots+\l_{j-1}+2,...,\l_1+\cdots+\l_{j-1}+\l_j\}$
for $1\le j\le k$ (we understand $\l_{-1}=0$). Let $S_\l$ be the
subgroup of $S_n$ consisting of elements stablizing each $I_j$. Then
$S_\l$ is a parabolic subgroup of $S_n$ and is isomorphic to
$S_{\l_1}\times S_{\l_2}\times\cdots\times S_{\l_k}$. We shall
denote by $w_\l$ the longest element of $S_\l$. Set $C_\l=\sum_{w\in
S_\l}T_w$. Following [KL] and [M] we consider the left idea
$N_\l=HC_\l$ of $H$ and shall regard it as a left $H$-module. Let
$N'_\l$ be the maximal submodule of $N_\l$ not containing $C_\l$.
Then the quotient module $M_\l=N_\l/N'_\l$ is an irreducible module
of $H$. Assume that $K$ is a field, then each irreducible module of
$H$ is isomorphic to some $M_\l$. See [KL, proof of Theorem 1.4] or
[M]. When $\sum_{w\in S_n}q^{l(w)}\ne 0$, the irreducible modules
$M_\l, \ \l$ a partition of $n$, form a complete set of irreducible
modules of $H$ (see [G, M], when $q$ is not a root of 1, this result
was implied in [L]).

One of the main result in [M] is the following.
\vskip3mm

\noindent{\bf Theorem.} Assume that $K$ is a field. Then

(a) The set $\{M_\l\ |\ C_\l M_\l\ne 0\}$ is a complete set of
irreducible modules of $H$.

(b) $C_\l M_\l\ne 0$ if and only if $\sum_{a=0}^mq^a\ne 0$ for all
$1\le m\le\text{max}\{\l_1-\l_2,\ \l_2-\l_3,\ ..., \l_{k-1}-\l_k,\
\l_k\}$. (See [M, Theorems 6.4 and 6.9]).

\vskip3mm

Now we argue for the theorem. For each module $E$ we can attached a
partition $\l=p(E)$ as follows, $C_\l E\ne 0$ but $C_\mu E=0$ for
all partition $\mu$ satisfying $\mu>\l$. (We say that
$\mu=(\mu_1,\mu_2,...,\mu_j)\ge\l=(\l_1,\l_2,...,\l_k)$ if
$\mu_1+\cdots+\mu_i\ge\l_1+\cdots+\l_i$ for $i=1,2,...$.)

Consider the two-sided ideal $F_\l=HC_\l H$ of $H$. According to the
proof of Theorem 1.4 in [KL], $F_\l/(F_\l\cap\sum_{\mu>\l}F_\mu)$ is
isomorphic to  the direct sum of some copies of
$E_\l=N_\l/(N_\l\cap\sum_{\mu>\l}F_\mu)$.

Let $E'_\l$ be the sum of all submodules $E$ of $E_\l$ satisfying
$C_\l E=0$. We claim that either $E'_\l=E_\l$ or $E'_\l$ is the
unique maximal submodule of $E_\l$.

 Let $D$ be a submodule of $E_\l$ such
that $C_\l D\ne 0$. For any $h\in H$ we have  $C_\l hC_\l\in aC_\l
+\sum_{\mu>\l}F_\mu$, here $a\in K$ ( loc.cit). Thus $C_\l D\ne 0$
implies that $C_\l D=E_\l$. Therefore $E'_\l=E_\l$ or $E'_\l$ is the
unique maximal submodule of $E_\l$. As a consequence,
$M_\l=E_\l/E'_\l$ if $C_\l E_\l\ne 0$ and in this case $C_\l M_\l\ne
0$.

Now assume that $L$ is an irreducible $H$-module such that $C_\l
L\ne 0$ but $C_\mu L=0$ for all $\mu>\l$.  Let $x\in L$ be such that
$C_\l x\ne 0$. Consider the $H$-module homomorphism $N_\l\to L,\
C_\l\to C_\l x$. By assumption,  $F_\mu L=0$ if $\mu>\l$. Thus we
get a nonzero homomorphism $E_\l\to L$. We must have $C_\l E_\l\ne
0$ since $C_\l L\ne 0$. So $L$ is isomorphic to $M_\l$. Noting that
$C_\mu E_\l\ne 0$ implies that $\mu\le\l$ (loc.cit) we see that  if
$\l\ne\mu$  then $M_\l$ is not isomorphic to $M_\mu$ when $C_\l
M_\l\ne 0\ne C_\mu M_\mu$. Part (a) is proved.

To prove part (b)   we need calculate $C_\l HC_\l$. This is
equivalent to calculate all  $C_\l T_wC_\l$. Clearly if $w\in S_\l$,
then $T_wC_\l=q^{l(w)}C_\l$. So we only need consider the element of
minimal length in a double coset $S_\l wS_\l$. Now the affine Hecke
algebra plays a role in calculating the product $C_\l T_wC_\l$.

Let $G$ be the special linear group $SL_n(\bold C)$ and let $T$ be
the subgroup of $G$ consisting of diagonal matrices. Let
$X=\text{Hom}(T,\bold C^*)$ be the character group of $T$. Let
$\tau_i\in X$ be the character $T\to \bold C,$
diag$(a_1,a_2,...,a_n)\to a_i$. Then we have
$\tau_1\tau_2\cdots\tau_n=1$ and as a free abelian group $X$ is
generated by $\tau_i,\ i=1,2,...,n-1$. The symmetric group $S_n$
acts on $X$ naturally: $w:\ X\to X,\ \tau_i\to \tau_{w(i)}$. Thus we
can form the semi-direct product $\tilde S_n=S\ltimes X$. In $\tilde
S_n$ we have $w\tau_{i}=\tau_{w(i)}w$ for $w$ in $S_n$. Let
$s_0=s\tau_1^2\tau_2\cdot\tau_i\cdot\tau_{n-1}$, where $s\in S_n$ is
the transpose $(1,n)=s_1s_2\cdots s_{n-2}s_{n-1}s_{n-2}\cdots
s_2s_1$. Since $\tau_1\tau_2\cdots\tau_n=1$ we have $s_0^2=1$. The
simple reflections $s_0,s_1,...,s_{n-1}$ generate a subgroup $W$ of
$\tilde S_{n}$, which is a Coxeter group of type $\tilde A_{n-1}$.
Define $\omega=\tau_1 s_1s_2\cdots s_{n-1}$. Then $\omega^n=1$ and
$\o s_i=s_{i+1}\o$ for all $i$ (we set $s_n=s_0$). Let $\Omega$ be
the subgroup of $\tilde S_n$ generated by $\o$. Note that $W$ is a
normal subgroup of $\tilde S_n$ and we have $\tilde
S_n=\Omega\ltimes W$. The Hecke algebra $\tilde H$ of $\tilde S_n$
is defined as follows. As a $K$-module, it is free and has a basis
consisting of elements $T_w,\ w\in \tilde S_n$. The multiplication
is defined by the relations $(T_{s_i}-q)(T_{s_i}+1)=0$ for all $i$
and $T_wT_u=T_{wu}$ if $l(wu)=l(w)+l(u)$. The length function $l:\
\tilde S_n\to \bold N$ is defined as $l(\omega^a w)=l(w)$ for $w\in
W$. Clearly $H$ is a subalgebra of $\tilde H$.

For $1\le i\le n-1$, define $x_i=\tau_1\tau_2\cdots\tau_i$. Then we
have $s_ix_j=x_js_i$ if $i$ and $j$ are different. Moreover we have
$l(w_0\prod_{i=1}^{n-1}x_i^{a_i})=l(w_0)+\sum_{i=1}^{n-1}a_il(x_i)$
if all $a_i$ are non-negative integers. Here $w_0$ is the longest
element of $S_n$. Also we have  $l(x_is_j)=l(x_i)-1$ if and only if
$i=j$.

Thus we have $T_{s_i}T_{x_j}=T_{x_j}T_{s_i}$ if $1\le i\ne j\le n-1$
and $T_{x_i}=T_{x_is_i}T_{s_{i}}$.

For a positive integer $k$ we set $[k]=q^{k-1}+q^{k-2}+\cdots+q+1$,
$[k]!=[k][k-1]\cdots[2][1]$, we also set $[0]=[0]!=1$.  For any
element $w\in \tilde S_n$ we set $C_w=\sum_{y\le w}P_{y,w}(q)T_y$,
where $\le$ is the Bruhat order and $P_{y,w}$ is the Kazhdan-Lusztig
polynomial. Note that if $w$ is a longest element of a parabolic
subgroup of $\tilde S_n$, then $C_w=\sum_{y\le w}T_y$. So we have
$C_\l=C_{w_\l}$. Now we are ready to prove part (b) of the theorem.

\vskip3mm

 \noindent{\bf Lemma 1.} Let $\l=(i,1,...,1)$ be a
partition of $n$ and $z\in S_n$ such that for  any simple reflection
$s$, $sz\le z$ if and only if $s=s_{i}$ and $zs\le z$ if and only if
$s=s_i$. Then
$$C_\l T_zC_\l\in \pm q^{*}[i-j-1]!C_\mu+\sum_\nu F_\nu,$$ for some
$j\le i-1$, where $*$ stands for an integer, $\mu=(i,j+1,1,...,1)$,
the summation runs through $\nu=(i+m,j+1-m,1,...,1)>\mu$ for $j+1\ge
m\ge 1$.

Proof: Since for any simple reflection $s$, if $sz\le z$ or $zs\le
z$ then we have $s=s_i$, we can find $j\le i-1$ such that
$$z=(s_is_{i-1}\cdots s_{i-j})(s_{i+1} s_i\cdots s_{i-j+1})\cdots
(s_{i+j-1}s_{i+j-2}\cdots s_{i-1})(s_{i+j}s_{i+j-1}\cdots s_{i}).$$
It is no harm to assume  $n=i+j+1$.

Note that
$$x_i=\omega^i (s_{n-i}s_{n-i-1}\cdots s_1)( s_{n-i+1}s_{n-i}\cdots
s_2)\cdots (s_{n-1}s_{n-2}\cdots s_i).$$

Let $y=(s_{i-j-1}s_{i-j}\cdots s_{i-1})\cdots (s_2s_3\cdots
s_{j+2})(s_1s_2\cdots s_{j+1})$. Since $n=i+j+1$ we have
$z=y\omega^{-i}x_i$ and $l(x_i)=l(y^{-1})+l(z)$ (we understand that
$y=e$ if $j=i-1$.) Thus we have $C_\l T_zC_\l=C_\l
T_{y^{-1}}^{-1}T_\omega^{-i}T_{x_i}C_\l.$ Noting that $C_\l
T_{y^{-1}}^{-1}=q^{-l(y)}C_\l$ and $C_\l T_{x_i}=T_{x_i}C_\l$, we
get
$$C_\l T_zC_\l=q^{-l(y)}C_\l
T_\omega^{-i}T_{x_i}C_\l=q^{-l(y)}T_\omega^{-i}T_\omega^{i}C_\l
T_\omega^{-i}C_\l T_{x_i}.$$ Let $w'=\omega^{i}w_\l\omega^{-i}$.
Then $w'$ is the longest element of the subgroup of $\tilde S_n$
generated by $s_{i+1},s_{i+2},...,s_{i+i-1}$. Let $k=i-j-2$, then
$2i-1=k+i+j+1$. We have $w'=u w_k$ for some $u$ and
$l(w')=l(u)+l(w_{k})$, where $w_{k}$ is the longest element of the
subgroup $W_k$ of $S_n$ generated by $s_1,s_2,...,s_k$ if $k\ge 1$
and $w_k=e$ is the neutral element if $k=-1$ or 0. We also have
$u=u'u_{i+1}$for some $u'$ and $l(u)=l(u')+l(u_{i+1})$, where
$u_{i+1}$ is the longest element of the subgroup of $U_{i+1}$ of
$S_n$ generated by $s_{i+1},...,s_{i+j}=s_{n-1}$. So
$T_\omega^{i}C_\l T_\omega^{-i}=hC_{u_{i+1}}C_{w_k}$ for some $h$ in
$H$, where $C_{u_{i+1}}$ is the sum of all $T_x$, $x\in U_{i+1}$,
and $C_{w_k}$ is the sum of all $T_x,\ x\in W_k$. Clearly we have
$C_{w_k}C_\l=[k+1]!C_\l$ and $C_{u_{i+1}}C_\l=C_\mu$. Therefore
$C'_\l C_\l=[k+1]!hC_{\mu}$. Note that $uw_\l=u'u_{i+1}w_\l=u'w_\mu$
is in the subgroup of $\tilde S_n$ generated by $s_p,\ p\ne i$. The
subgroup is isomorphic to the symmetric group $S_n$. Applying the
Robinson-Schensted rule we see that $uw_\l$ and $w_\mu$ are in the
same left cell. (See [A] for an exposition of Robinson-Schensted
rule. One may see this fact also from star operations introduced in
[KL].) Write $C'_\l C_\l=\sum a_vC_v$, then clearly
$a_{uw_\l}=[k+1]!$. Since $T_\omega$ and $T_{x_{i}}$ are invertible,
we see that in the expression $C_\l T_zC_\l=\sum b_vC_v,\ b_v\in K$,
there exists $x$ such that $b_x\ne 0$, $x$ and $w_\mu$ are in the
same two-sided cell. Since $z=z^{-1}$ and $w_{\l}=w_\l^{-1}$, by the
symmetry we see that $x$ and $w_\mu$ are in the same left cell  and
right cell as well. So we must have $x=w_\mu$ (see [KL, proof of
Theorem 1.4]). Moreover we must have $b_\mu=\pm q^a[k+1]!$ for some
integer $a$. If $b_v\ne 0$ and $v\ne w_{\mu}$, we must have $C_v\in
F_\nu$ for some $\nu>\mu$. We claim that for such $\nu$ we have
$\nu=(i+m,j+1-m,1,...,1)$  for some $m\ge 1$ . Since $C_\l T_zC_\l$
is contained in the subalgebra of $H$ generated by
$T_{s_1},...,T_{s_{i+j}}$, we may assume that $n=i+j+1$. In this
case we must have $\nu=(i+m,j+1-m)$ for some $m\ge 1$ since
$\mu=(i,j+1)$ and $\nu>\mu$. The lemma is proved.

\vskip3mm

Remark: The author has not been able to determine the integer $*$ in
the lemma.

\vskip3mm

\noindent{\bf Corollary 2.} Let $\l=(i,j)$ be a partition of $n$.
That is $i\ge j$ and $i+j=n$. Then for any $z$ in $S_n$ we have
$$C_\l T_zC_\l\in [i-j]![j]!fC_\l+\sum_{\mu>\l}F_{\mu},$$
where $f\in K$.

Proof: Since $C_\l C_\l=[i]![j]!C_\l$ and $T_sC_\l=C_\l T_s=qC_\l$
if $s\ne i$ in $S_n$, we may assume that $z= (s_is_{i-1}\cdots
s_{i-k})\cdots (s_{i+k}s_{i+k-1}\cdots s_i)$, where $k\le j-1\le
i-1$. Note that $C_\l=C_{w_{i-1}}C_{u_{i+1}}=C_{u_{i+1}}C_{w_{i-1}}$
(see the proof of Lemma 1 for the definition of $w_i$ and $u_i$). We
have $C_\l T_zC_\l
=C_{u_{i+1}}C_{w_{i-1}}T_zC_{w_{i-1}}C_{u_{i+1}}$. By Lemma 1 we get
$C_{w_{i-1}}T_zC_{w_{i-1}}\in\pm
q^*[i-k-1]!C_{w_{i-1}w_{i+1,i+k}}+\sum_\nu F_{\nu}$, where
$w_{i+1,i+k}$ is the longest element of the subgroup of
$W_{i+1,i+k}$ of $S_n$ generated by $s_{i+1},...,s_{i+k}$, and $\nu$
runs through the partitions $(i+m,k+1-m,1,...,1),\ k+1\ge m\ge 1$.

We have
$C_{u_{i+1}}C_{w_{i+1,i+k}}C_{u_{i+1}}=[k+1]![j]!C_{u_{i+1}}$. We
also have $C_\l T_zC_\l \subset\sum_{\mu\ge\l}F_{\mu}$ and
$C_{u_{i+1}}F_{\nu}C_{u_{i+1}}\subset \sum_{\mu\ge \nu} F_{\mu}$ for
any $\nu$. If $\mu\ge\l$ and $\mu\ge (i+m,...,)$ for some $m\ge 1$,
we must have $\mu>\l$. So $C_\l T_zC_\l\in \pm
q^*[i-k-1]![k+1]![j]!C_\l+\sum_{\mu>\l}F_{\mu}.$ Since
$[i-k-1]!=[i-j]![i-j+1]\cdots[i-k-1]$, the corollary follows.

\vskip3mm

\noindent{\bf Lemma 3.} Let $\l=(\l_1,\l_2,...,\l_k)$ be a partition
of $n$. Then
$$C_\l T_zC_\l\in\prod_{i=1}^k[\l_i-\l_{i+1}]!fC_\l+F_{>\l},$$ where $f\in K$ and we set $\l_{k+1}=0.$

Proof: We use induction on $k$. When $k=1$, the lemma is trivial,
when $k=2$, by  Corollary 2 we see  the assertion  is true. Now
assume that $k>2$. For $i\le j$ we set $l_{i,j}=\l_i+\cdots \l_j$.
We have (see the proof of Corollary 2 for the definition of
$w_{km}$)
$$w_\l=w_{\l_1-1}w_{\l_1+1,\l_{1,2}-1}\cdots
w_{\l_{1,k-1}+1,\l_{1,k}-1}=w_{\l_1-1}w'.$$ Let $z=xz_1y$, where
$x,y$ are in the subgroup of $S_n$ generated by $s_i,\ i\ne \l_1$,
and $l(s_iz_1)=l(z_1s_i)=l(z_1)+1$ if $i\ne \l_1$. Write $x=x_1x_2$
and $y=y_1y_2$, where $x_1,y_1$ are in the subgroup $W_{\l_1-1}$ of
$S_n$ generated by $s_1,...,s_{\l_1-1}$ and $x_2,y_2$ are in the
subgroup $U_{\l_1+1}$ of $S_n$ generated by
$s_{\l_1+1},...,s_{n-1}$.

 We  have $T_{u}C_\l=C_\l
T_{u}=q^{l(u)}=q^{l(u)}C_\l$ for $u=x_1,y_1$ and
$T_{u}C_{w_{\l_1-1}}=C_{w_{\l_1-1}}T_u$ for $u=x_2,y_2$. Note that
$C_\l=C_{w_{\l_1-1}}C_{w'}=C_{w'}C_{w_{\l_1-1}}$. Thus
$$C_\l T_zC_\l=q^{l(x_1)+l(y_1)}C_{w'}T_{x_2}C_{w_{\l_1-1}}T_{z_1}C_{w_{\l_1-1}}T_{y_2}C_{w'}.$$
If $z_1=e$, then $$C_\l
T_zC_\l=q^{l(x_1)+l(y_1)}[\l_1]!C_{w_{\l_1-1}}C_{w'}T_{x_2y_2}C_{w'}.$$
We are reduced to the case $k-1$.

Now assume that $z_1\ne e$. By Lemma 1 we know that
$$C_{w_{\l_1-1}}T_{z_1}C_{w_{\l_1-1}}\in \pm
q^*[\l_1-j-1]C_{w_{\l_1-1}w_{\l_1+1,\l_1+j}}+\sum_\nu F_{\nu},$$
where $j\le\l_1-1$ is defined by $z_1= s_{\l_1}s_{\l_1-1}\cdots
s_{\l_1-j}\cdots s_{\l_1+j}s_{\l_1+j-1}\cdots s_{\l_1},$  and $\nu$
runs through the partitions $(\l_1+m,j+1-m,1,...,1),\ j+1\ge m\ge
1$.

Note that both $C_\l T_zC_\l$ and
$C_{w'}T_{x_2}C_{w_{\l_1-1}w_{\l_1+1,\l_1+j}}T_{y_2}C_{w'} $ are
contained in $\sum_{\mu\ge\l}F_{\mu}$ and
$C_{w'}T_{x_2}F_{\nu}T_{y_2}C_{w'}\subset \sum_{\mu\ge \nu}
F_{\tau}$ for any $\nu$. Whenever $\mu\ge\l$ and
 $\mu\ge (\l_1+m,...,)$ for some $m\ge
1$, we must have $\mu>\l$. Thus we have $$C_\l T_zC_\l\in \pm
q^*[\l_1-j-1]!C_{w'}T_{x_2}C_{w_{\l_1-1}w_{\l_1+1,\l_1+j}}T_{y_2}C_{w'}+\sum_{\mu>\l}F_\mu,$$
where $*$ stands for an integer.  Let $\tau=(\l_1,j+1,1,...,1).$
Then $w_{\l_1-1}w_{\l_1+1,\l+j}=w_\tau.$ Note that
$C_{w_{\l_1-1}w_{\l_1+1,\l_1+j}}=C_{\l_1-1}C_{w_{\l_1+1,\l_1+j}}.$
If $j\ge \l_2$, then $\tau\not\le\l$, so
$C_{w'}T_{x_2}C_{w_{\l_1-1}w_{\l_1+1,\l_1+j}}T_{y_2}C_{w'}$ is
contained in
$(\sum_{\mu\ge\l}F_\mu)\cap(\sum_{\mu\ge\tau}F_\mu)\subset\sum_{\mu>\l}F_\mu.$
We are done in this case.

Now assume that $j\le\l_2-1$. Then $\l_1-j-1\ge\l_1-\l_2$ and
$C_{w_{\l_1-1}}T_{z_1}C_{w_{\l_1-1}}\in
[\l_1-\l_2]!f_1C_{w_{\l_1-1}w_{\l_1+1,\l+j}}+\sum_{\nu}F_{\nu}$ for
some $f_1\in K$, where $\nu$ runs through the partitions
$(\l_1+m,j+1-m,1,...,1)$, $j+1\ge m\ge 1$.  Thus we have
$$C_{\l}T_zC_\l\in
[\l_1-\l_2]!f_1C_{w_{\l_1-1}}C_{w'}T_{x_2}C_{w_{\l_1+1,\l+j}}T_{y_2}C_{w'}
+\sum_{\nu}C_{w'}T_{x_2}F_{\nu}T_{y_2}C_{w'}.$$ Note that $x_2, w',
y_2, w_{\l_1+1,\l+j}$ are all in the subgroup of $S_n$ generated by
$s_i,\ \l_1+1\le i\le n-1$ and $C_{w'}T_{x_2}F_{\nu}T_{y_2}C_{w'}$
is included in $\sum_{\mu\not\le\l}F_\mu$ if
$\nu=(\l_1+m,j+1-m,1,...,1)$ for some $m\ge 1$. By induction
hypothesis, we see the lemma is true.

\vskip3mm

\noindent{\bf Lemma 4.} Let $\l$ be as in Lemma 3. Set $$z_i=
(s_{\l_{1i}}s_{\l_{1i}-1}\cdots s_{\l_{1i}-\l_{i+1}+1})\cdots
(s_{\l_{1,i+1}-1}s_{\l_{1,i+1}-2}\cdots s_{\l_{1,i}}),$$ for
$i=1,2,...,k-1$. Define
$$h=T_{z_{k-1}}(T_{z_{k-2}}T_{z_{k-1}})(T_{z_{k-3}}T_{z_{k-2}}T_{z_{k-1}})\cdots
(T_{z_{1}}T_{z_{2}}\cdots T_{z_{k-1}}).$$ Then $C_\l hC_\l\in\pm
q^*\prod_{i=1}^k([\l_i-\l_{i+1}]!)^iC_\l+F_{>\l},$ where $*$ stands
for an integer and $F_\l=\sum_{\mu>\l}F_\mu$.

Proof: Set $u_i=C_{w_{\l_{1,i-1}+1,\l_{1i}-1}}$ (we understand that
$\l_{1,0}=0$) and $h_i=T_{z_i}$. Then $C_\l=u_1u_2\cdots u_k$,
$u_iu_j=u_ju_i$ for all $i,j$, and $u_ih_j=h_ju_i$ if $i<j$. For
$h',h''\in H$ and $F\subset H$, we write $h'\equiv h''+F$ if
$h'-h''\in F$. Using Lemma 1 we get
 $$\begin{array}{rl} C_\l hC_\l &=u_k(u_{k-1}h_{k-1})(u_{k-2}h_{k-2}h_{k-1})\cdots
\\[3mm] &\quad \times (u_2h_2h_3\cdots h_{k-1})u_1h_1u_1h_2u_2\cdots h_{k-1}u_{k-1}u_k\\[3mm]
&\equiv\pm q^* [\l_1-\l_2]!u_k(u_{k-1}h_{k-1})(u_{k-2}h_{k-2}h_{k-1})\cdots\\[3mm] &
\quad \times (u_2h_2h_3\cdots h_{k-1})u_1u_2h_2u_2\cdots h_{k-1}u_{k-1}u_k+F_{>\l}\\[3mm]
&\equiv\pm q^* [\l_1-\l_2]!u_1u_k(u_{k-1}h_{k-1})(u_{k-2}h_{k-2}h_{k-1})\cdots\\[3mm] &
\quad \times (u_2h_2h_3\cdots h_{k-1})u_2h_2u_2\cdots h_{k-1}u_{k-1}u_k+F_{>\l}\\[3mm]
&\equiv\pm q^*
[\l_1-\l_2]![\l_2-\l_3]!u_1u_k(u_{k-1}h_{k-1})(u_{k-2}h_{k-2}h_{k-1})\cdots\\[3mm]
& \quad\times
(u_2h_2h_3\cdots h_{k-1})u_2u_3h_3u_3\cdots h_{k-1}u_{k-1}u_k+F_{>\l}\\[3mm]
&\equiv\pm q^*
[\l_1-\l_2]!([\l_2-\l_3]!)^2u_1u_2u_k(u_{k-1}h_{k-1})(u_{k-2}h_{k-2}h_{k-1})\cdots
\\[3mm]&\quad\times (u_3h_3\cdots h_{k-1})^2u_3h_3u_3\cdots
h_{k-1}u_{k-1}u_k+F_{>\l}\\[3mm]
&\equiv\cdots\\[3mm]
&\equiv\pm q^*\prod_{i=1}^k([\l_i-\l_{i+1}]!)^iC_\l+F_{>\l}.
\end{array}$$

Combining Lemmas 3 and 4 we see that part (b) of the theorem is
true. The theorem is proved.

\vskip3mm

If $\sum_{w\in S_n}q^{l(w)}\ne 0$ and $K$ is an algebraic closed
field of characteristic 0, then we have the
Deligne-Langlands-Lusztig classification for irreducible modules of
$\tilde H$ (see [BZ, Z], [KL1], [X]). We have another classification
due to Ariki and Mathas for any sufficient large $K$ (see [AM]). An
interesting question is  to classify irreducible modules of $\tilde
H$ in the spirit of Deligne-Langlands-Lusztig classification when
$\sum_{w\in S_n}q^{l(n)}= 0$, see [Gr] for an announcement. If one
can manage the calculation $C_\l\tilde HC_\l$ to get counterparts of
Lemmas 3 and 4 , the question will be settled.

\vskip3mm

\noindent{\bf Acknowledgement:} Part of the paper was written during
my visit to the  Department of Mathematics at the National
University of Singapore. I  am grateful to Professors C. Zhu and D.
Zhang for invitation and   to the department for hospitality and
financial support.

\end{document}